\documentclass[12pt]{article}
\usepackage{amssymb,latexsym, amsmath, amscd, array, graphicx}
\usepackage{amssymb}

\usepackage{setspace}

\usepackage{circledsteps}

\input xy
\xyoption{all}
\makeatletter 

\usepackage{setspace}
\usepackage{enumitem}

\date{}

\newtheorem{theorem}{Theorem}[section]

\newtheorem{proposition}[theorem]{Proposition}

\newcommand{\z}{{\Bbb Z}}
\newcommand{\D}{{\Bbb D}}

\newcommand{\re}{{\Bbb R}}

\newcommand{\lo}{\rightarrow}

\newcommand{\black}{{\blacksquare}}

\newcommand{\diam}{{\rm diam}}

\newcommand{\mdim}{{\rm mdim}}
\newcommand{\mesh}{{\rm mesh}}
\newcommand{\ord}{{\rm ord}}

\begin{document}

\title{\bf An embedding theorem  for mean dimension}

\author{  Michael  Levin
\footnote{The author
 was supported by the Israel Science Foundation grant No. 2196/20 and 
 the Fields institute of mathematics   where 
 this note  was written  during his visit in 2023.}}

\maketitle

\begin{abstract}
Let   $(X,\z)$  be a minimal dynamical system   on a compact metric $X$  and $k$ an integer such that $\mdim X< k$.
 We show that $(X,\z)$ admits an equivariant embedding in the shift $(\D^k)^\z$ where $\D$ is a superdendrite.
\\
{\bf Math. Subj. Class.:}  37B05 (54F45)
\end{abstract}
\begin{section}{Introduction}
This note is devoted to proving 

\begin{theorem}
\label{theorem-dendrite}
 Let $(X, \z)$ be a minimal dynamical system on a compact metric $X$
  and let $k$ be a natural number such that $\mdim X < k$. Then  for  almost every map $f : X \lo \D^k$
 the induced map  $f^\z : X \lo  (\D^k)^\z$ is an embedding where $\D$ is a superdendrite.

\end{theorem}
We recall that a superdendrite  is a dendrite  with  a dense set of end-points \cite{bower,sternfeld}.
A superdendrite $\D$ is a $1$-dimensional compact metric AR (absolute retract), $\D$ is embeddable in the plane and
$\D^\z$ is homeomorphic to the Hilbert cube \cite{bower,sternfeld,west}.

Since $\D$ is $1$-dimensional, 
$\mdim (\D^k)^\z =k$ \cite{tsukamoto}. Since $\D$ is embeddable in the plane, 
$(\D^k)^\z$ is equivariantly embeddable in $([0,1]^{2k})^\z$ and therefore
the Lindenstrauss-Tsukamoto examples  \cite{lindenstrauss-tsukamoto} show 
that  the inequality $\mdim X < k$ in Theorem \ref{theorem-dendrite} cannot be improved to the sharp
inequality $\mdim X \leq k$.

This note is based on the approach of \cite{levin-1} and we adopt the notations of \cite{levin-1}.

\end{section}

\begin{section}{Preliminaries}
We  will present here some facts and notations  used in the proof of Theorem \ref{theorem-dendrite}.
\begin{theorem}
\label{theorem-sternfeld}
{\rm (\cite{sternfeld}, Theorem 2.1)}
Let $X$ be compact metric and $a \in X$.
Then for almost every map $f : X \lo \D$  to a superdendrite $\D$ we have $f^{-1}(f(a))=\{ a\}$.
\end{theorem}

\begin{proposition}

\label{prop-dendrite}
Let $f : X \lo \D$ be a map from a compact metric space $X$ to a superdendrite $\D$. Then for every $\alpha>0$
there is $\beta>0$ such that for every finite collection $\cal A$ of closed disjoint subsets of $X$ with $\mesh \cal A < \beta$
the map $f$ can be approximated by an $\alpha$-close map $g: X \lo \D$ sending the elements of $\cal A$ to distinct singletons in $\D$
such that $g^{-1}(g(A))=A$ for every $A \in \cal A$.
\end{proposition}
{\bf Proof.} Since $\D$ is a compact metric AR there is $\beta>0$ such that for every surjective  map $\phi : X \lo Y$ to a compact metric $Y$
such that  the fibers of $\phi$ are of $\diam < \beta$  there is a map $\psi : Y\lo \D$ such that $\psi\circ \phi$ is $\alpha/2$-close to $f$.
We will  show that the proposition  holds for this $\beta$. Let $Y$ be obtained from $X$ by collapsing the elements of $\cal A$ to
singletons and let $\phi : X \lo Y$ be the projection. Then there is a map $\psi : Y \lo \D$ such that $\psi\circ\phi$ is $\alpha/2$-close to $f$.
By Theorem \ref{theorem-sternfeld}, we can replace $\psi$ by an $\alpha/2$-close map and 
assume that $\psi\circ \phi$ is $\alpha$-close to $f$ and $\psi^{-1}(\psi(a))=\{a\}$ for every $\{a\}=f(A), A\in \cal A$.
Set $g=\psi \circ \phi$ and the proposition follows.
$\black$\\

Let  $(Y, \re)$ be a dynamical system, $A$ a subset of $Y$,  $\cal A$ a collection of subsets of $Y$ and $\alpha, \beta \in \re$ positive numbers.
The subset $A$ is said to be {\bf $(\alpha, \beta)$-small if $\diam( A +r) < \alpha$} for every $r \in [-\beta, \beta]\subset \re$.
The collection  $\cal A$ is said to be {\bf $(\alpha, \beta)$-fine} if  $\mesh({\cal A}+r)< \alpha$ for every $r \in [0, \beta]\subset  \re$.
The collection $\cal A$ is said to be {\bf   $(\alpha, \beta)$-refined} at  a subset $W \subset Y$ if the following two conditions hold:
(condition 1) no element  of ${\cal A}+r$ meets
the closures of both  $W+r_1$ and $W+r_2$ for every $r,r_1, r_2 \in [-\beta, \beta]\subset \re$  with $|r_1-r_2|\geq 1$ and 
(condition 2) if for an element $A$ of ${\cal A}$ the set $A+[-\beta,\beta]$ meets the closure of $W+[-\beta, \beta]$ 
then $\diam (A+r) < \alpha$ for every $r \in [-\beta, \beta]\subset \re$.

\begin{proposition}
\label{prop-1}
{\rm (\cite{levin-1})}
Let $(Y, \re)$ be a free dynamical system on a compact metric $Y$, $\omega$ a point in $Y$ and let $\alpha$ and  $\beta $ be positive real numbers.
Then there is an open neighborhood  $W$ of $w$  and an open cover $\cal V$ of $Y$ such that $\ord {\cal V} \leq 3$ and $\cal V$ is $(\alpha, \beta)$-refined at $W$. 
\end{proposition}

\begin{proposition}
\label{prop-2} 
{\rm (\cite{levin-1})}
Let $q>2$ be an integer. Then there is  a finite collection ${\cal E}$ of  disjoint closed intervals in $ [0, q) \subset \re$
such that $\cal E$  splits into the union ${\cal  E}={\cal E}_1 \cup \dots \cup {\cal E}_q$  of $q$ disjoint subcollections 
having the property that for every $ t \in \re$  the set  $t +\z\subset \re$  meets at least  $q-2$   subcollections ${\cal E}_i$
(a set meets a collection  if there is a point of the set that is covered by the collection). Moreover,  we may assume
that $\mesh {\cal E}$ is as small as we wish. 
\end{proposition}

\begin{theorem}
\rm{ (\cite{levin-1})}
\label{theorem-borel}
For any dynamical system $(X,\z)$ on a compact metric $X$  one has  \\
$\mdim X \times_\z \re=\mdim X$ where $X \times_\z \re$  is  Borel's construction for $(X,\z)$.
\end{theorem}

We recall that Borel's construction   $X \times_\z \re$ is also known as the mapping torus  in topological dynamical. 
The space $X$ is  naturally  embedded in $X \times_\z \re$ so that the action of $\z$ in $X \times_\z \re$
extends the action of $\z$ on $X$, and  Borel's construction  $X \times_\z  \re$ is endowed with the standard $\re$-action
that extends the $\z$-action  on $X \times_\z \re$ \cite{levin-1}.

\end{section}
\begin{section}{Proof of Theorem \ref{theorem-dendrite}}
Let   $f=(f_1, \dots, f_k):X \lo \D^k$ be any  map.  Fix $\delta>0$.  Apply Proposition \ref{prop-dendrite} with $\alpha=\delta$
to produce $\beta$ and set $\epsilon =\beta/3$.
Our  goal is to approximate $f$ by a $\delta$-close map  $\psi$ 
such that  the  fibers of $\psi^\z$ are of $\diam <3\epsilon$. 

By Theorem \ref{theorem-borel} we have $\mdim X \times_\z \re =\mdim X $. Let $q>3$ be a natural number and
$n=(q-3)k$. Recall $\mdim X < k$. Then, assuming that $q$ is large enough, there is 
an open cover ${\cal U}$ of $X \times_\z \re$  such that $\ord {\cal U} \leq n-2$ and $\cal U$ is $(\epsilon, q)$-fine.

Since the theorem obviously  holds if $X$ is a singleton, we may assume that $(X,\z)$ is  non-trivial. Fix a point $w \in X $.
By Proposition \ref{prop-1} there is an open cover $\cal V$ of $X\times_\z \re$ and a neighborhood $W$ of $w$ in $X \times_\z \re$ such
 that $\cal V$ is $(\epsilon, 3q)$-refined at $W$. 
 
 Now replacing 
 $\cal U$ by an open cover of $\ord\leq n$  refining ${\cal U}\vee {\cal V}$ 
 we can assume that ${\ord {\cal U}} \leq n$, $\cal U$ is $(\epsilon, q)$-fine
and  $\cal U$ is $(\epsilon,3q)$-refined at $W$. Clearly we can replace $W$ by a smaller neighborhood of $w$ and assume
that $W$ is $(\epsilon, 3q)$-small and  the elements of ${\cal D}_W$ are disjoint where
${\cal D}_W$ is the collection of the closures of $W+z$ for the integers 
 $z\in [-3q, 3q]$.
 
Set $m=qk$. Refine $\cal U$ by a  Kolmogorov-Ostrand cover
  ${\cal F}$ of $X \times_\z \re$   such that  
 ${\cal F}$ covers $X \times_\z \re$ at least $m-n=3k$ times and
 ${\cal F}$ splits into ${\cal F}={\cal F}_1 \cup \dots \cup {\cal F}_{m}$ the union of finite  families of
 disjoint closed sets ${\cal F}_i$. Note that $\cal F$ is $(\epsilon,q)$-fine and $(\epsilon, 3q)$-refined at $W$.

 Let $\xi : X \lo \re$ be a Lindenstrauss level function determined by $W $ restricted to $X$.
Denote $W^+ =W +\z \cap [-q,q]$ and $X^-=X \setminus W^+$.
Recall  that 
$\xi(x+z)=\xi(x) +z$ for every  $x \in X^-$ and  an integer $-q\leq z\leq q$.

Following \cite{levin-1} we need an additional auxiliary notation.
Let  $\cal A$ be  a collection of subsets of $X \times_\z \re$, $\cal B$ a collection of intervals in $\re$. 
For $B\in \cal B$ and $z \in \z$ consider the collection ${\cal A}+B$ restricted to $\xi^{-1}(B+qz)$ and
denote by ${\cal A}\oplus_\xi B$ the union of such collections for all $z\in \z$.   Now  denote
by ${\cal A}\oplus_\xi {\cal B}$ the union
of the collections ${\cal A}\oplus_\xi  B$ for all $ B \in {\cal B}$. Note that ${\cal A}\oplus_\xi {\cal B}$ is
a collection of subsets of $X$.

Consider  a finite  collection ${\cal E}$ of  disjoint closed intervals in $ [0, q) \subset \re$ satisfying  the conclusions of Proposition \ref{prop-2}.
 For $1\leq i \leq k$  define the collection  ${\cal D}_i$  of subsets of $X$ as the union of the collections
${\cal F}_{i}\oplus_\xi {\cal E}_1$, ${\cal F}_{i+k}\oplus_\xi {\cal E}_2$, \dots, ${\cal F}_{i+(q-1)k}\oplus_\xi {\cal E}_q$. 
Note that assuming that $\mesh \cal E$ is small enough we may also assume that ${\cal F}_i^+={\cal F}_i+[-\mesh{\cal E}, \mesh{\cal E}]$
is a collection of disjoint sets and 
 the collection ${\cal F}^+={\cal F}+[-\mesh{\cal E}, \mesh{\cal E}]$ is 
$(\epsilon,q)$-fine and $(\epsilon, 3q)$-refined at $W$  and, as a result,  we get that
 ${\cal D}_i $ is a collection of disjoint  closed sets of $X$ of $\diam < \epsilon$ and
each  element of ${\cal D}_i$ meets at most one element  of ${\cal D}_W$.
 
Then, by Proposition \ref{prop-dendrite}, we can define 
a map $\psi=(\psi_1, \dots, \psi_k) : X\lo \D^k$  so that  for each $i$ the map  $\psi_i$ is $\delta$-close to $f_i$,
$\psi_i$ sends the elements of ${\cal D}_W$ restricted to $X$  and the elements of  ${\cal D}_i$  to singletons in $\D$, the  preimage
of each such singleton under $\psi_i$ is exactly the union of the elements of ${\cal D}_W$ and ${\cal D}_i$ sent by $\psi_i$ to that singleton and, finally,
$\psi_i$ separates  the elements of ${\cal D}_W$ restricted to $X$  together with  the elements of ${ \cal D}_i$ not  meeting ${\cal D}_W$. 
We will show that the fibers of $\psi^\z$ are of  $\diam <  3\epsilon$.

Consider  a point $x \in X\setminus X^-$. Then $x$ is contained in an element $D$ of ${\cal D}_W$ and
therefore $\diam \psi_i ^{-1}(\psi_i (x))< 3\epsilon$ for every $i$. 

Now consider a point $x \in X^-$ and let $z, l\in \z$ be such that $zq \leq l\leq  \xi(x) <l +1\leq z(q+1)$.
Set $y=x-(l-zq)$ and note that $zq\leq \xi(y)< zq+1$.
Recall  that  the point  $y+(zq-\xi(y))\in X \times_\z \re$ is covered by at least $3k$ collections 
from the family ${\cal F}_1, \dots, {\cal F}_m$ and $\xi(y)+\z$ meets   at least 
$q-2$ collections from ${\cal E}_1, \dots, {\cal E}_q$.  Then there is an integer $i$  and a collection ${\cal D}_j$ such that
$0\leq i < q$  and  the point  $y +i$ is covered by ${\cal D}_j$.

Indeed, let  $\xi(y) +\z$  meet ${\cal E}_p$. Then, since
$zq\leq \xi(y)<zq +1$, there is an integer $i_p$ such that  $0\leq i_p <q$, $zq \leq \xi(y) +i_p <z(q+1)$ and $\xi(y)+i_p$ is covered by
${\cal E}_p+zq$.  
 Note that   different $p$ define different $i_p$ and
for every $1\leq j\leq k$ such that ${\cal F}_{j+(p-1)k}$ covers the point  $y+(zq-\xi(y))$ we have 
that the collection ${\cal D}_j$ covers $y+i_p$.
Thus if $\xi(y)+\z$ meets all the collections ${\cal E}_1, \dots, {\cal E}_q$  
the number of   collections ${\cal D}_j$ meeting $y+i$ for some integer  $0\leq i <q$
will be 
at least the number of times $y +(zq-\xi(y))$ is covered 
by the collections ${\cal F}_1, \dots, {\cal F}_m$,  which is
at least $3k$.
Each time  $\xi(y)+\z$ misses a collection from ${\cal E}_1, \dots, {\cal E}_q$  reduces the above estimate by at most $k$. 
Since $\xi(y)+\z$ can miss at most two collections from ${\cal E}_1, \dots, {\cal E}_q$ 
we get there is an integer $0\leq i <q$ and  a collection ${\cal D}_j$ such that ${\cal D}_j$ 
covers $y +i$
and $zq\leq \xi(y)+i<(z+1)q$. 

Let  $D \in {\cal D}_j$ be the element containing $y+i$.
Note that $D$ is contained in an element of 
${\cal F}^+  + \xi(y)-zq +i$.

Assume that $D$ meets  ${\cal D}_W$. Then $\diam\psi_j^{-1}(\psi_j(D)) <3\epsilon$. Moreover,
since ${\cal D}_W$ is $(\epsilon, 3q)$-small and 
and ${\cal F}^+$ is $(\epsilon,  3q)$-refined at $W$ we get that 
$\diam (\psi_j^{-1}(\psi_j (D))+t) <3\epsilon$ for every real $-3q\leq t \leq 3q$. Thus
 we have $x\in D +(l-zq) -i$ and  $(\psi_j^\z)^{-1}(\psi_j^\z(x)) \subset \psi_j^{-1}(\psi_j (D))+(l-zq)-i$ and
 hence   $\diam (\psi_j^\z)^{-1}(\psi_j^\z(x))<3\epsilon$ since $-3q\leq (l-zq)-i\leq 3q$.

Now assume that $D$ does not meet   ${\cal D}_W$. Then $\psi_j^{-1}(\psi_j(D)) =D$.
Recall that  $D$ is contained in an element of ${\cal F}^+ + \xi(y)-zq +i={\cal F}^+ +\xi(x-(l-zq)) -zq +i ={\cal F}^+ +\xi(x) -l +i$.
  Then $D +(l-zq)-i$ is contained in an element of ${\cal F}^+ +\xi(x)-zq$ and note  that  $x \in D +(l-zq)-i$.
Since ${\cal F}^+$ is $(\epsilon,q)$-fine and $ 0 \leq \xi(x)-zq < q$ we get $\diam  (D +(l-zq)-i) < \epsilon$.
Then, since $(\psi_j^\z)^{-1}(\psi_j^\z(x)) \subset \psi_j^{-1}(\psi_j (D))+(l-zq)-i=D+(l-zq)-i$, we get
$\diam (\psi_j^\z)^{-1}(\psi_j^\z(x))<\epsilon$.

Thus for every $x \in X$ there is $j$ such that $\diam (\psi_j^\z)^{-1}(\psi_j^\z(x))<3\epsilon$ and hence 
the fibers of $\psi^\z$ are of $\diam < 3\epsilon$ and the theorem follows by a standard Baire  category argument.
$\black$

\end{section}

 Department of Mathematics\\
Ben Gurion University of the Negev\\
P.O.B. 653\\
Be'er Sheva 84105, ISRAEL  \\
 mlevine@math.bgu.ac.il\\\\ 
  
\end{document}